# A NEW MAXIMAL INEQUALITY AND INVARIANCE PRINCIPLE FOR STATIONARY SEQUENCES

By Magda Peligrad[1] and Sergey Utev

*University of Cincinnati and University of Nottingham*

We derive a new maximal inequality for stationary sequences under a martingale-type condition introduced by Maxwell and Woodroofe [*Ann. Probab.* **28** (2000) 713–724]. Then, we apply it to establish the Donsker invariance principle for this class of stationary sequences. A Markov chain example is given in order to show the optimality of the conditions imposed.

**1. Results.** Let $(X_i)_{i \in Z}$ be a stationary sequence of centered random variables with finite second moment ($E[X_1^2] < \infty$ and $E[X_1] = 0$). Denote by $\mathcal{F}_k$ the $\sigma$-field generated by $X_i$ with indices $i \le k$, and define

$$S_n = \sum_{i=1}^{n} X_i, \qquad W_n(t) = \frac{S_{[nt]}}{\sqrt{n}}, \qquad 0 \le t \le 1,$$

where $[x]$ denotes the integer part of $x$. Finally, let $W = \{W(t) : 0 \le t \le 1\}$ be a standard Brownian motion. In the sequel $\Longrightarrow$ denotes the weak convergence and $\|X\| = \sqrt{E(X^2)}$.

Theorem 1.1. *Assume that*

$$(1) \qquad \sum_{n=1}^{\infty} \frac{\|E(S_n|\mathcal{F}_0)\|}{n^{3/2}} < \infty.$$

*Then, $\{\max_{1 \le k \le n} S_k^2/n : n \ge 1\}$ is uniformly integrable and $W_n(t) \Longrightarrow \sqrt{\eta} W(t)$, where $\eta$ is a nonnegative random variable with finite mean $E[\eta] = \sigma^2$ and*

Received May 2003; revised April 2004.

[1]Supported in part by a Taft grant.

*independent of* $\{W(t); t \geq 0\}$. *Moreover, condition* (1) *allows to identify the variable* $\eta$ *from the existence of the following limit*

$$\lim_{n \to \infty} \frac{E(S_n^2|\mathcal{I})}{n} = \eta \qquad in \ L_1,$$

*where* $\mathcal{I}$ *is the invariant sigma field. In particular,* $\lim_{n \to \infty} E(S_n^2)/n = \sigma^2$.

In the next theorem we show that, in its generality, condition (1) is optimal in the following sense.

THEOREM 1.2. *For any nonnegative sequence* $a_n \to 0$, *there exists a stationary ergodic discrete Markov chain* $(Y_k)_{k \geq 0}$ *and a functional* $g$ *such that* $X_i = g(Y_i)$; $i \geq 0$, $E[X_1] = 0$, $E[X_1^2] < \infty$ *and*

$$\sum_{n=1}^{\infty} a_n \frac{\|E(S_n|Y_0)\|}{n^{3/2}} < \infty \qquad but \ \frac{S_n}{\sqrt{n}} \ is \ not \ stochastically \ bounded.$$

In the ergodic case, Theorem 1.1 improves upon the corresponding results of Maxwell and Woodroofe (2000) [see also Derriennic and Lin (2003) and Wu and Woodroofe (2002)].

Our method of proof is based on the martingale approximation originated in Gordin (1969). Rather then considering and analyzing a perturbed solution of the Poisson equation, as it was suggested in Maxwell and Woodroofe (2000) [see also Liverani (1996)], we analyze small blocks and apply maximal inequalities to show that the sums of variables in these blocks can be approximated by stationary martingale differences.

In the proof of our key inequalities, we use a variety of techniques. The starting point is the diadic induction found to be useful in the analysis of $\rho$-mixing sequences. This method goes back to Ibragimov (1975), and was further developed by many authors including Peligrad (1982), Shao (1989), Bradley and Utev (1994) and Peligrad and Utev (1997). The second tool is the modification of the Garsia (1965) telescoping sums approach to maximal inequalities as used by Peligrad (1999) and Dedecker and Rio (2000). Our maximal inequality, stated in Proposition 2.3, is new and has interest in itself. Finally, we use the subadditivity of the conditional sums of random variables.

In order to show the optimality of our results, we construct an example which is motivated by the well-known counterexample stating that, in the general ergodic case, unlike the i.i.d. case (the Kolmogorov strong law of the large numbers), $E|X| = \infty$ does not imply that the averages $S_n/n$ diverge almost surely [see Halmos (1956), page 32; he has attributed this example to M. Gerstenhaber]. The discrete version of the example was probably introduced in Chung [(1960), Markov chains, page 92]. For the modern development and connection with Pomeau–Manneville type 1 intermittency model, we mention Isola (1999) whose detailed analysis was inspirational.

Theorem 1.1 is proved in Sections 2.1–2.4. Theorem 1.2 is proved in Sections 3.1–3.3.

**2. Proof of Theorem 1.** Throughout the section we will use the notation

$$\Delta_r = \sum_{j=0}^{r-1} \left\| \frac{E(S_{2^j}|\mathcal{F}_0)}{2^{j/2}} \right\|. \tag{4}$$

### 2.1. *Analysis of second-order moments of partial sums.*

PROPOSITION 2.1. *Let $n$, $r$ be integers such that $2^{r-1} < n \leq 2^r$. Then*

$$E(S_n^2) \leq n[\|X_1\| + \tfrac{1}{2}\Delta_r]^2. \tag{5}$$

*Assume $\sum_{j=0}^{\infty} 2^{-j/2}\|E(S_{2^j}|\mathcal{F}_0)\| < \infty$. Then, the following limit exists in $L_1$:*

$$\eta := \lim_{n \to \infty} \frac{E(S_n^2|\mathcal{I})}{n} = E(X_1^2|\mathcal{I}) + \sum_{j=0}^{\infty} \frac{E[S_{2^j}(S_{2^{j+1}} - S_{2^j}|\mathcal{I})]}{2^j}, \tag{6}$$

*where $\mathcal{I}$ is the invariant sigma field. In particular,*

$$\sigma^2 := E[\eta] = E(X_1^2) + \sum_{j=0}^{\infty} \frac{E(S_{2^j}(S_{2^{j+1}} - S_{2^j}))}{2^j}.$$

PROOF. The last statement is an immediate consequence of (6). In order to prove (5), we shall use an induction argument. It is easy to see that (5) is true for $r = 0$ and $n = 1$. Assume (5) holds for all $n \leq 2^{r-1}$. Fix $n$, $2^{r-1} < n \leq 2^r$. Starting with $S_n = S_{n-2^{r-1}} + S_n - S_{n-2^{r-1}}$ and using the Cauchy–Schwarz inequality and stationarity, we derive

$$\|S_n\|^2 \leq \|S_{n-2^{r-1}}\|^2 + \|S_{2^{r-1}}\|^2 + 2\|S_{n-2^{r-1}}\|\|E(S_{2^{r-1}}|\mathcal{F}_0)\|.$$

Now, by induction assumption, since $\|E(S_{2^{r-1}}|\mathcal{F}_0)\| = 2^{(r-1)/2}(\Delta_r - \Delta_{r-1})$, and $4(n - 2^{r-1})2^{r-1} \leq n^2$, we obtain

$$\|S_n\|^2 \leq (n - 2^{r-1})[\|X_1\| + \tfrac{1}{2}\Delta_{r-1}]^2 + 2^{r-1}[\|X_1\| + \tfrac{1}{2}\Delta_{r-1}]^2$$
$$+ 2(n - 2^{r-1})^{1/2}[\|X_1\| + \tfrac{1}{2}\Delta_{r-1}]2^{(r-1)/2}(\Delta_r - \Delta_{r-1})$$
$$\leq n[\|X_1\| + \tfrac{1}{2}\Delta_{r-1} + \tfrac{1}{2}(\Delta_r - \Delta_{r-1})]^2 = n[\|X_1\| + \tfrac{1}{2}\Delta_r]^2.$$

This establishes (5).

To prove (6) for the subsequence $n = 2^r$, we use the notation $E_I(Y) = E(Y|\mathcal{I})$ and $\|Y\|_I = \sqrt{E(Y^2|\mathcal{I})}$ for the corresponding norm. By recurrence, we can easily establish the representation

$$
\begin{aligned}
E_I(S_{2^r}^2) &= 2^r E_I(X_1^2) + \sum_{i=1}^{r} 2^i E_I[S_{2^{r-i}}(S_{2^{r-i+1}} - S_{2^{r-i}})] \\
&= 2^r \left( E_I(X_1^2) + \sum_{j=0}^{r-1} \frac{E_I(S_{2^j}(S_{2^{j+1}} - S_{2^j}))}{2^j} \right).
\end{aligned}
$$
(7)

We observe that

$$
E[S_{2^j}(S_{2^{j+1}} - S_{2^j})|\mathcal{I}] = E\{E[S_{2^j}(S_{2^{j+1}} - S_{2^j})|\mathcal{F}_{2^j}]|\mathcal{I}\}
$$

[see, e.g., Proposition (2.2) in Bradley (2002), page 54]. Thus, by the Jensen inequality,

$$
E|E[S_{2^j}(S_{2^{j+1}} - S_{2^j})|\mathcal{I}]| \le E|E[S_{2^j}(S_{2^{j+1}} - S_{2^j})|\mathcal{F}_{2^j}]|
$$

so that the Cauchy–Schwarz inequality and stationarity imply

$$
E|E[S_{2^j}(S_{2^{j+1}} - S_{2^j})|\mathcal{I}]| \le \|S_{2^j}\|\|E(S_{2^j}|\mathcal{F}_0)\|.
$$

In addition, by the first part of the proposition and the summability of the series in the right-hand side of (4), we obtain

$$
\sum_{j=0}^{\infty} \frac{\|S_{2^j}\|\|E(S_{2^j}|\mathcal{F}_0)\|}{2^j} \le C \sum_{j=0}^{\infty} \left\| \frac{E(S_{2^j}|\mathcal{F}_0)}{2^{j/2}} \right\| < \infty,
$$

which proves the convergence in $L_1$ of the series

$$
E(X_1^2|\mathcal{I}) + \sum_{j=0}^{\infty} \frac{E[S_{2^j}(S_{2^{j+1}} - S_{2^j})|\mathcal{I}]}{2^j} = \eta.
$$

This relation and (7) show that the convergence in (6) holds along the subsequence $n = 2^r$, that is,

$$
\lim_{r \to \infty} \frac{E[S_{2^r}^2|\mathcal{I}]}{2^r} = \eta.
$$

To treat the whole sequence $S_n$, for $1 \le n < 2^r$, we start with the binary expansion

$$
n = \sum_{k=0}^{r-1} 2^k a_k \qquad \text{where } a_{r-1} = 1 \text{ and } a_k \in \{0, 1\}.
$$

Then, we apply the following representation:

$$
S_n = \sum_{j=0}^{r-1} T_{2^j} a_j \qquad \text{where } T_{2^j} = \sum_{i=n_{j-1}+1}^{n_j} X_i, n_j = \sum_{k=0}^{j} 2^k a_k, n_{-1} = 0.
$$

Clearly, for $a_j = 0$, $T_{2^j} = 0$. For $a_j = 1$, the conditional distribution of $T_{2^j}$ given $\mathcal{I}$ is equally distributed as the conditional distribution of $S_{2^j}$ given $\mathcal{I}$.

To prove (6), we start with the representation

$$E[S_n^2 | \mathcal{I}] = \left( \sum_{i=1}^{r-1} a_i E[S_{2^i}^2 | \mathcal{I}] \right) + \left( \sum_{\substack{i \neq j=1}}^{r-1} a_i a_j E[T_{2^i} T_{2^j} | \mathcal{I}] \right) \equiv I_n + J_n.$$

By the above convergence, $E[S_{2^j}^2 | \mathcal{I}]/2^j \to \eta$, which implies the convergence

$$\frac{I_n}{n} \to \eta \qquad \text{in } L_1.$$

It remains to prove that $\frac{E|J_n|}{n} \to 0$. Let $i < j < r$. Then, as before,

$$E|E[T_{2^i} T_{2^j} | \mathcal{I}]| \leq E|E[T_{2^i} E(T_{2^j} | \mathcal{F}_{n_i})]| \leq \|S_{2^i}\| \|E(S_{2^j} | \mathcal{F}_0)\|$$

$$\leq C 2^{i/2} \sqrt{n} \left\| \frac{E(S_{2^j} | \mathcal{F}_0)}{2^{j/2}} \right\|$$

and, thus,

$$E|J_n| \leq 2 \sum_{1=i<j\leq r-1} E|E[T_{2^i} T_{2^j} | \mathcal{I}]| \leq 2C\sqrt{n} \sum_{i=1}^{r-2} 2^{i/2} \sum_{j=i+1}^{r} \left\| \frac{E(S_{2^j} | \mathcal{F}_0)}{2^{j/2}} \right\|,$$

which implies $E|J_n|/n \to 0$ because

$$\sum_{j=i}^{\infty} \left\| \frac{E(S_{2^j} | \mathcal{F}_0)}{2^{j/2}} \right\| \to 0 \qquad \text{as } i \to \infty. \qquad \square$$

2.2. *Maximal inequalities.* We start by establishing first an auxiliary lemma:

LEMMA 2.2. *Let $(Y_i)_{1 \leq i \leq n}$ be a random vector of square integrable random variables such that for each $i$, $1 \leq i \leq n$, $Y_i$ is measurable with respect to $\mathcal{F}_i = \sigma(X_j, j \leq i)$, where $(X_i)$ is a stationary sequence introduced before. Let $n \leq 2^r$. If for all $1 \leq a \leq b \leq n$, and a positive constant $C$,*

$$E \left( \sum_{l=a}^{b} Y_l \right)^2 \leq C(b - a + 1) \qquad then \quad \left| E \sum_{l=1}^{n-1} Y_l (S_n - S_l) \right| \leq \tfrac{1}{2} C n \Delta_r.$$

PROOF. We shall prove this lemma by induction. It is easy to see the result of this lemma is true for $n = 2$. Assume the lemma holds for all $n \leq 2^{r-1}$. Fix now $n$, $2^{r-1} < n \leq 2^r$, and begin by writing

$$\sum_{l=1}^{n-1} Y_l (S_n - S_l) = \sum_{l=1}^{n-2^{r-1}-1} Y_l (S_{n-2^{r-1}} - S_l)$$

$$+ \sum_{l=n-2^{r-1}}^{n-1} Y_l(S_n - S_l) + \sum_{l=1}^{n-2^{r-1}-1} Y_l(S_n - S_{n-2^{r-1}}).$$

$$= I_1 + I_2 + I_3.$$

By using the Cauchy–Schwarz inequality, along with the conditions of this lemma and stationarity, we easily obtain

$$|EI_3| \leq C[2^{r-1}(n - 2^{r-1})]^{1/2}(\Delta_r - \Delta_{r-1}) \leq \tfrac{1}{2}Cn[\Delta_r - \Delta_{r-1}].$$

By the induction assumption, $|EI_1| \leq \tfrac{1}{2}C(n - 2^{r-1})\Delta_{r-1}$ and $|EI_2| \leq \tfrac{1}{2}C2^{r-1} \times \Delta_{r-1}$, so

$$|EI_1| + |EI_2| + |EI_3| \leq \tfrac{1}{2}Cn\Delta_{r-1} + \tfrac{1}{2}Cn[\Delta_r - \Delta_{r-1}] = \tfrac{1}{2}Cn\Delta_r,$$

proving the lemma.   □

We are ready to state and prove our key maximal inequality.

PROPOSITION 2.3.   *Let* $\{X_i : i \in Z\}$ *be a stationary sequence of random variables. Let* $n$, $r$ *be integers such that* $2^{r-1} < n \leq 2^r$. *Then we have*

$$E\left[\max_{1 \leq i \leq n} S_i^2\right] \leq n(2\|X_1\| + (1 + \sqrt{2})\Delta_r)^2.$$

PROOF.   Denote by

$$M_n = \max_{1 \leq i \leq n} |S_i| \quad \text{and} \quad K_m = \max_{1 \leq j \leq m} \frac{1}{j} E\left[\max_{1 \leq i \leq j} S_i^2\right].$$

We first prove that, for any positive integer $n$,

(8)        $$E\left[\max_{1 \leq i \leq n} S_i^2\right] \leq n[2K_n^{1/2}\Delta_r + 4[\|X_1\| + \tfrac{1}{2}\Delta_r]^2].$$

By the fact that $K_l$ is nondecreasing in $l$, from (8), we derive

$$K_n \leq 2K_n^{1/2}\Delta_r + 4[\|X_1\| + \tfrac{1}{2}\Delta_r]^2,$$

which implies $K_n^{1/2} \leq 2\|X_1\| + (1 + \sqrt{2})\Delta_r$, hence, the result.

To prove (8), we denote by $S_0 = 0$,

$$M_n^+ = \max_{1 \leq j \leq n} S_j^+ = \max(0, S_1, \ldots, S_n)$$

and

$$M_n^- = \max_{1 \leq j \leq n} (-S_j^-) = \max(0, -S_1, \ldots, -S_n).$$

We shall use the following simplified version of an interesting inequality in Dedecker and Rio (2000) [see (3.4) in Dedecker and Rio (2000) or (3.5) in

Rio (2000)], which was obtained by using Garsia's (1965) telescoping sum approach to the maximal inequality

$$(9) \qquad (M_n^+)^2 \leq 4(S_n^+)^2 - 4 \sum_{k=1}^{n} M_{k-1}^+ X_k.$$

By adding to this relation the similar one for $M_n^-$, we obtain

$$(M_n)^2 \leq 4(S_n)^2 - 4 \sum_{k=1}^{n} (M_{k-1}^+ - M_{k-1}^-)(X_k).$$

We now write $X_k = (S_n - S_{k-1}) - (S_n - S_k)$ and derive

$$(10) \qquad (M_n)^2 \leq 4(S_n)^2 - 4 \sum_{k=1}^{n-1} D_k(S_n - S_k),$$

where $D_k = (M_k^+ - M_{k-1}^+) - (M_k^- - M_{k-1}^-)$.

It is easy to see that

$$\left| \sum_{k=a+1}^{b} D_k \right| \leq \max[(M_b^+ - M_a^+), (M_b^- - M_a^-)]$$

$$\leq \max_{a \leq i \leq b} |S_i - S_a|.$$

Taking the expectation, we get, by stationarity,

$$E\left( \sum_{k=a+1}^{b} D_k \right)^2 \leq E\left( \max_{1 \leq i \leq b-a} S_i^2 \right) = (b-a)K_{b-a} \leq (b-a)K_n.$$

Next, by Lemma 2.2 applied with $Y_k = D_k$ for $k \geq 1$, $C = K_n^{1/2}$, we obtain

$$\left| E \sum_{k=1}^{n-1} D_k(S_n - S_k) \right| \leq \tfrac{1}{2} n[K_n^{1/2} \Delta_r].$$

By substituting this estimate in (10) together with (5) on $E(S_n^2)$, we obtain (8) and, hence, the proposition. □

REMARK 2.4. The inequality in Proposition 2.3 is an extension of the Doob maximal inequality for martingales, giving also an alternative proof of this famous theorem. Notice that, for the martingale case, our inequality gives the same constant as in the Doob inequality, a constant that cannot be improved. A natural question that arises is the optimality of the constant in front of $\Delta_r$ and further study is needed to determine the best constants in this inequality.

2.3. *Analysis of certain series involving conditional sums.*

(a) *Key result.* Let $X = (X_i)_{i \in Z}$ be a stationary sequence of random variables with finite second moment. Denote by

$$S_n = \sum_{i=1}^{n} X_i, \qquad V_n = V_n(X) = \|E(S_n|\mathcal{F}_0)\|,$$

where as before, $\mathcal{F}_k$ is the $\sigma$-field generated by $X_i$ with indices $i \leq k$.

The main condition (1) of Theorem 1.1 is $\sum V_n/n^{3/2} < \infty$. On the other hand, various inequalities derived in Sections 2.1 and 2.2 have used the condition $\sum V_{2^r}/2^{r/2} < \infty$. In this section we show that these conditions are equivalent and, in addition, we prove the following proposition, which is useful in establishing the martingale approximation in Theorem 1.1.

PROPOSITION 2.5. *Under condition* (1),

$$\frac{\|E(S_m|\mathcal{F}_0)\|}{\sqrt{m}} \to 0 \quad \text{and} \quad \frac{1}{\sqrt{m}} \sum_{j=0}^{\infty} \left\| \frac{E(S_{m2^j}|\mathcal{F}_0)}{2^{j/2}} \right\| \to 0$$

*as* $m \to \infty$.

PROOF. In order to prove this result, we shall analyze in Lemma 2.6 the conditional variance of sums and then, in Lemma 2.7, some related series. By Lemma 2.6, the sequence $V_m = \|E(S_m|\mathcal{F}_0)\|$ is subadditive. Then, we have only to apply Lemma 2.8 to conclude the proof of this proposition. □

(b) *Conditional variances of sums form a subadditive sequence.* The starting point of our analysis is the following simple observation.

LEMMA 2.6. $V_n$ *is a subadditive sequence.*

PROOF. First, since for all $n$, $\mathcal{F}_{-n} \subset \mathcal{F}_0$, we observe that

$$E[E(S_k|\mathcal{F}_{-n})]^2 \leq E[E(S_k|\mathcal{F}_0)]^2 = \|E(S_k|\mathcal{F}_0)\|^2 = V_k^2.$$

Hence, by stationarity,

$$\|E(S_{i+j} - S_i|\mathcal{F}_0)\| = \sqrt{E[E(S_j|\mathcal{F}_{-i})]^2} \leq V_j.$$

Thus,

$$V_{i+j} = \|E(S_i + [S_{i+j} - S_i]|\mathcal{F}_0)\| \leq \|E(S_i|\mathcal{F}_0)\| + \|E(S_{i+j} - S_i|\mathcal{F}_0)\|$$
$$\leq V_i + V_j.$$

(c) *Analysis of certain series for subadditive sequences.* Let $V_n$ be a nonnegative subadditive sequence. For a $p > 1$, define

$$I := \sum_{j=0}^{\infty} \frac{V_{2^j}}{2^{j(p-1)}}, \qquad J := \sum_{n=1}^{\infty} \frac{V_n}{n^p}, \qquad W := \sum_{n=1}^{\infty} n^{-p} \max_{1 \leq i \leq n} V_i.$$

$\square$

The following lemma is a crucial step in deriving the result in Proposition 2.5.

LEMMA 2.7. *There exists two positive absolute constants $C_p$ and $K_p$ such that*

$$C_p I \leq J \leq W \leq K_p I.$$

PROOF. We shall start with the following simple representation:

$$W = \sum_{n=1}^{\infty} n^{-p} \max_{1 \leq i \leq n} V_i = \sum_{r=0}^{\infty} \sum_{n=2^r}^{2^{r+1}-1} n^{-p} \max_{1 \leq i \leq n} V_i.$$

Then, by the subadditivity of the sequence $\{V_n; n \geq 0\}$, for $i \leq n < 2^{r+1}$,

$$V_i \leq \sum_{j=0}^{r} V_{2^j} \qquad \text{so that} \max_{1 \leq i \leq n} V_i \leq \sum_{j=0}^{r} V_{2^j},$$

which implies

$$W \leq \sum_{r=0}^{\infty} 2^{-pr} 2^r \sum_{k=0}^{r} V_{2^k} = \sum_{k=0}^{\infty} V_{2^k} \sum_{r=k}^{\infty} 2^{-r(p-1)} = K_p \sum_{k=0}^{\infty} 2^{-k(p-1)} V_{2^k}$$

$$= K_p I,$$

where $K_p = \frac{1}{1-2^{-(p-1)}}$. The last inequality is therefore proved.

The inequality $J \leq W$ is straightforward. Now, we need the following simple combinatorial property. Define

$$A_N = \{1 \leq i \leq N : V_i \geq V_N/2\} \qquad \text{and denote by } |A| \text{ the cardinal of a set } A.$$

PROPERTY. $|A_N| \geq N/2$, *that is, $A_N$ contains at least $N/2$ elements.*

PROOF. To prove it, we denote by $D_N = \{1, \ldots, N\}$ and fix $1 \leq i < N$. Observe that if $i \in A_N^c = D_N - A_N$, then $N - i \in A_N$ because

$$V_{N-i} \geq V_N - V_i > V_N - V_N/2 \geq V_N/2.$$

Thus, $A_N \supseteq N - A_N^c$ and so $N = |D_N| = |A_N| + |A_N^c| \leq 2|A_N|$ and the property is proved. $\square$

Now, in order to continue the proof of Lemma 2.7, we write

$$J = \sum_{r=0}^{\infty} \left( \sum_{n=4^r}^{4^{r+1}-1} \frac{V_n}{n^p} \right) \geq \frac{1}{4^p} \sum_{r=0}^{\infty} 4^{-rp} \left( \sum_{n=4^r}^{4^{r+1}-1} V_n \right).$$

We are going to apply the above property with $N = 4^{r+1}$. Define

$$C_r = \{n \in \{4^r, \dots, 4^{r+1} - 1\} : V_n \geq V_N/2\} = A_N \cap \{4^r, \dots, 4^{r+1} - 1\}.$$

Clearly,

$$|C_r| \geq |\{4^r, \dots, 4^{r+1} - 1\}| - |A_N^c| = 4^{r+1} - 4^r - |A_N^c|$$

and, applying the above property, we obtain

$$|C_r| \geq 4^{r+1} - 4^r - (4^{r+1} - 1)/2 \geq 4^{r+1} - 4^r - 4^{r+1}/2 = 4^r.$$

Thus,

$$J \geq \frac{1}{2} \frac{1}{4^p} \sum_{r=0}^{\infty} 4^{-rp} V_{4^{r+1}} |C_r| \geq \frac{1}{2} \frac{1}{4^p} \sum_{r=0}^{\infty} 4^{-r(p-1)} V_{4^{r+1}} = \frac{1}{8} \sum_{r=1}^{\infty} 2^{-2r(p-1)} V_{2^{2r}},$$

which implies

$$Q := \sum_{r=0}^{\infty} 2^{-2r(p-1)} V_{2^{2r}} = V_1 + \sum_{r=1}^{\infty} 2^{-2r(p-1)} V_{2^{2r}} \leq 9J.$$

Then, by subadditivity, $V_{2^{2r+1}} \leq 2 V_{2^{2r}}$, so that

$$P := \sum_{r=0}^{\infty} 2^{-(2r+1)(p-1)} V_{2^{2r+1}} \leq \frac{2}{2^{(p-1)}} \sum_{r=0}^{\infty} 2^{-(2r)(p-1)} V_{2^{2r}} = \frac{2}{2^{(p-1)}} Q$$

and, as a consequence,

$$I = \sum_{r=0}^{\infty} \frac{V_{2^{2r}}}{2^{2r(p-1)}} + \sum_{r=0}^{\infty} \frac{V_{2^{2r+1}}}{2^{(2r+1)(p-1)}} = P + Q \leq 9\left(\frac{2}{2^{(p-1)}} + 1\right) J,$$

and the proof of Lemma 2.7 is complete.   $\square$

LEMMA 2.8.   *Assume that* $\sum_{n=1}^{\infty} V_n n^{-3/2} < \infty$. *Then,*

$$G_m = \frac{1}{\sqrt{m}} \sum_{k=0}^{\infty} \frac{V_{m2^k}}{2^{k/2}} \to 0 \qquad \text{as } m \to \infty.$$

*In particular,* $V_m/\sqrt{m} \to 0$ *as* $m \to \infty$.

PROOF.   By rewriting $G_m$, we obtain

$$G_m = \sum_{k=0}^{\infty} \sum_{n=m2^k}^{m2^{k+1}-1} (m2^k)^{-3/2} V_{m2^k} \leq 2^{3/2} \sum_{k=0}^{\infty} \sum_{n=m2^k}^{m2^{k+1}-1} n^{-3/2} \max_{1 \leq i \leq n} V_i$$

$$= 2^{3/2} \sum_{n=m}^{\infty} n^{-3/2} \max_{1 \leq i \leq n} V_i,$$

which proves that $G_m \to 0$ as $m \to \infty$ by Lemma 2.7.   $\square$

2.4. *Martingale approximation and the proof of Theorem* 1. Let $m$ be a fixed integer and $k = [n/m]$, where, as before, $[x]$ denotes the integer part of $x$. We start the proof by dividing the variables in blocks of size $m$ and making the sums in each block

$$X_i^{(m)} = m^{-1/2} \sum_{j=(i-1)m+1}^{im} X_j, \qquad i \geq 1.$$

Then we construct the martingale

$$M_k^{(m)} = \sum_{i=1}^{k} (X_i^{(m)} - E(X_i^{(m)} | \mathcal{F}_{i-1}^{(m)})), \qquad i \in Z,$$

where $\mathcal{F}_k^{(m)}$ denotes the $\sigma$-field generated by $X_i^{(m)}$ with indices $i \leq k$.

Notice that $M_k^{(m)}$ is a stationary martingale and, therefore, by the classical invariance principle for martingales, we derive

$$\frac{1}{\sqrt{k}} M_{[kt]}^{(m)} \Longrightarrow \sqrt{\eta^{(m)}} W,$$

where $\eta^{(m)}$ is the following limit (both in $L_1$ and almost surely):

$$\eta^{(m)} = \lim_{k \to \infty} \frac{1}{k} \sum_{i=1}^{k} (X_i^{(m)} - E(X_i^{(m)} | \mathcal{F}_{i-1}^{(m)}))^2.$$

In order to prove the invariance principle for $\frac{1}{\sqrt{n}} S_{[nt]}$, together with the uniform integrability of the sequence $\max_{1 \leq k \leq n} S_k^2/n$, by the Doob maximal inequality and Theorem 4.2 in Billingsley (1968), we have only to establish that

(11) $$\|\sqrt{\eta^{(m)}} - \sqrt{\eta}\| \to 0 \qquad \text{as } m \to \infty$$

and

(12) $$\lim_{m \to \infty} \lim_{n \to \infty} \left\| \sup_{0 \leq t \leq 1} \left| \frac{1}{\sqrt{n}} S_{[nt]} - \frac{1}{\sqrt{k}} M_{[kt]}^{(m)} \right| \right\| = 0.$$

Notice first that by the convergence in Proposition 2.5,

$$\lim_{m \to \infty} \frac{1}{m} E[E(S_m | \mathcal{F}_0)]^2 = 0.$$

On the other hand, by the ergodic theorem (both almost surely and in $L_1$),

$$\lim_{k \to \infty} \frac{1}{k} \sum_{i=1}^{k} (X_i^{(m)})^2 = \frac{1}{m} \lim_{n \to \infty} \frac{1}{n} \sum_{i=1}^{n} (S_{i+m} - S_i)^2 = \frac{E[S_m^2 | \mathcal{I}]}{m},$$

where $\mathcal{I}$ is the $\sigma$-field of invariant sets.

Therefore, by Proposition 2.1, we obtain the following convergence in $L_1$:

$$\lim_{m \to \infty} \eta^{(m)} = \lim_{m \to \infty} \frac{E(S_m^2|\mathcal{I})}{m} = \eta,$$

which implies (11).

To prove (12), we first notice that

$$\left\| \sup_{0 \le t \le 1} \left| \frac{1}{\sqrt{n}} S_{[nt]} - \frac{1}{\sqrt{km}} S_{[nt]} \right| \right\| \le \left( 1 - \frac{\sqrt{n}}{\sqrt{km}} \right) \left\| \frac{1}{\sqrt{n}} \max_{1 \le j \le n} |S_j| \right\|.$$

By taking into account Proposition 2.3 and the fact that $\lim_{n \to \infty}(1 - \frac{\sqrt{n}}{\sqrt{km}}) = 0$, the right-hand side of the above inequality tends to 0. Therefore, we have only to estimate

$$\left\| \sup_{0 \le t \le 1} \left| \frac{1}{\sqrt{km}} S_{[nt]} - \frac{1}{\sqrt{k}} M_{[kt]}^{(m)} \right| \right\|$$

$$\le \frac{1}{\sqrt{km}} \left\| \sup_{0 \le t \le 1} \sum_{i=[kt]m+1}^{[nt]} X_i \right\| + \frac{1}{\sqrt{k}} \left\| \sup_{0 \le t \le 1} \left| \sum_{i=1}^{[kt]} E(X_i^{(m)}|\mathcal{F}_{i-1}^{(m)}) \right| \right\|,$$

which leads to the estimate

$$\left\| \sup_{0 \le t \le 1} \left| \frac{1}{\sqrt{km}} S_{[nt]} - \frac{1}{\sqrt{k}} M_{[kt]}^{(m)} \right| \right\|$$

$$\le \frac{3m}{\sqrt{km}} \left\| \max_{1 \le i \le n} X_i \right\| + \frac{1}{\sqrt{k}} \left\| \max_{1 \le j \le k} \left| \sum_{i=1}^{j} E(X_i^{(m)}|\mathcal{F}_{i-1}^{(m)}) \right| \right\|.$$

Since for every $\epsilon > 0$,

$$E \max_{1 \le i \le n} X_i^2 \le \epsilon + \sum_{i=1}^{n} X_i^2 I(|X_i| > \epsilon)$$

by stationarity, for any fix $m$, $\lim_{n \to \infty} 3m \| \max_{1 \le i \le n} X_i \| / \sqrt{km} = 0$.

On the other hand, by Propositions 2.3 and 2.5, we derive

$$\frac{1}{\sqrt{k}} \left\| \max_{1 \le j \le k} \left| \sum_{i=1}^{j} E(X_i^{(m)}|\mathcal{F}_{i-1}^{(m)}) \right| \right\|$$

$$\le 2 \frac{\|E(S_m|\mathcal{F}_0)\|}{\sqrt{m}} + (1 + \sqrt{2}) \frac{1}{\sqrt{m}} \sum_{j=0}^{\infty} \left\| \frac{E(S_{m2^j}|\mathcal{F}_0)}{2^{j/2}} \right\| \to 0$$

as $m \to \infty$, uniformly in $n$, which completes the proof of Theorem 1.1.

## 3. Proof of Theorem 2.

3.1. *The countable Markov chain and its preliminary analysis.* Let $\{Y_k; k \geq 0\}$ be a discrete Markov chain with the state space $Z^+$ and transition matrix $P = (p_{ij})$ given by $p_{k(k-1)} = 1$ for $k \geq 1$ and $p_j = p_{0(j-1)} = P(\tau = j)$, $j = 1, 2, \ldots$ (i.e., whenever the chain hits $0$, $Y_t = 0$, it then regenerates with the probability $p_j$). When $p_1, p_2 > 0$, and, in addition, $p_{n_j} > 0$ along $n_j \to \infty$, the chain is irreducible and aperiodic. The stationary distribution exists if and only if $E[\tau] < \infty$ and it is given by

$$\pi_j = \pi_0 \sum_{i=j+1}^{\infty} p_i, \qquad j = 1, 2, \ldots,$$

where $\pi_0 = 1/E[\tau]$.

Let us consider now an arbitrary nonnegative sequence $a_n \to 0$ as in our Theorem 1.2. Notice that, without loss of generality, it is enough to assume that $a_n$ is a strictly decreasing sequence of real positive numbers.

The choice of $p_j$ further depends on this arbitrary nonnegative sequence $a_n$. First, we define a sequence $\{u_k; k = 1, 2, \ldots\}$ of positive integers such that

$$(13) \qquad u_1 = 1, \qquad u_2 = 2, \qquad u_k^4 + 1 < u_{k+1} \qquad \text{for } k \geq 3 \quad \text{and}$$

$$a_t \leq k^{-2} \qquad \text{for } t \geq u_k.$$

Then, for $i \geq 1$, we take

$$p_i = \begin{cases} c/u_j^2, & \text{if } i = u_j \text{ for some } j \geq 1, \\ 0, & \text{if } i \neq u_j \text{ for all } j \geq 1, \end{cases}$$

that is, for each positive integer $j \geq 1$, $p_{u_j} = c/u_j^2$ and $p_i = 0$ for $u_j < i < u_{j+1}$.

Clearly,

$$(14) \qquad\qquad E[\tau] < \infty \qquad \text{but } E[\tau^2] = \infty.$$

As a functional $g$, we take $I_{(x=0)} - \pi_0$, where $\pi_0 = P_\pi(Y_0 = 0)$ under the stationary distribution denoted by $P_\pi$ ($E_\pi$ denotes the expectations for the process started with the stationary distribution). The stationary sequence is defined by

$$X_j = I_{(Y_j=0)} - \pi_0 \qquad \text{so that } S_n = \sum_{j=1}^{n} X_j = \sum_{j=1}^{n} I_{(Y_j=0)} - n\pi_0.$$

By $P_k$ and $E_k$, we denote the probability and the expectation operator when the Markov chain is started at $k$, that is, $P(Y_0 = k) = 1$. Let

$$\nu = \min\{m \geq 1 : Y_m = 0\}, \qquad A_n = E_0[S_n], \qquad x \wedge y = \min(x, y).$$

PROPOSITION 3.1.

$$V_n = \|E(S_n|Y_0)\| \le \|\nu \wedge n\| + \max_{1 \le i \le n} |A_i|$$

$$\equiv I_n + J_n,$$

*where* $\|x\|^2 = \sum_{k=0}^{\infty} x_k^2 \pi_k$.

PROOF.   We first notice that $|S_n| \le n$ and $P_k(\nu = k) = 1$, so that, conditionally on $Y_0 = k$ (with $0 < k \le n$),

$$E_k(S_n) = E_k(S_k) + E_k(S_n - S_k).$$

The first term is bounded by $k$ and the second term is equal to $E_0(S_{n-k+1})$ since $Y_k = 0$. Thus,

$$|E_k(S_n)| \le k \wedge n + |A_{n-k+1}|. \qquad \square$$

3.2. *Proving that* $\sum a_n \|E(S_n|Y_0)\| n^{-3/2} < \infty$. By Proposition 3.1, it is enough to prove that

$$(15) \qquad \sum_{n=1}^{\infty} a_n I_n/n^{3/2} + \sum_{n=1}^{\infty} a_n J_n/n^{3/2} < \infty.$$

The first sum is easily treated by a straightforward analysis. Indeed, to analyze $I = \sum a_n I_n/n^{3/2}$, we first notice that, for $u_{t-1} \le j$,

$$\pi_j = \pi_0 \sum_{i=j+1}^{\infty} p_i \le \pi_0 c_1/u_t^2.$$

Therefore, we write, for $u_k < n \le u_{k+1}$,

$$I_n^2 = E_\pi(\nu \wedge n)^2 = \sum_{j=1}^{n} j^2 \pi_j + n^2 \sum_{j=n+1}^{\infty} \pi_j$$

$$\le \left[ \sum_{t=1}^{k} \left( \sum_{j=u_{t-1}+1}^{u_t} j^2 \pi_j \right) \right] + \left( \sum_{j=u_k+1}^{n} j^2 \pi_j \right) + n^2 \sum_{t=k+1}^{\infty} \left( \sum_{j=u_{t-1}+1}^{u_t} \pi_j \right)$$

$$\le c_2 \left[ \sum_{t=1}^{k} u_t^{-2} \left( \sum_{j=1}^{u_t} j^2 \right) \right] + \frac{c_3}{u_{k+1}} \left( \sum_{j=u_k+1}^{n} j^2 \right) + c_3 n^2 \sum_{t=k+1}^{\infty} \frac{1}{u_t}$$

$$\le c_4 (u_k + n^3/u_{k+1}^2 + n^2/u_{k+1}).$$

Next, write

$$\sum_{n=1}^{\infty} \frac{a_n I_n}{n^{3/2}} = \sum_{k=1}^{\infty} \sum_{n=u_k+1}^{u_{k+1}} \frac{a_n I_n}{n^{3/2}} \le \sum_{k=1}^{\infty} \frac{1}{k^2} \sum_{n=u_k+1}^{u_{k+1}} I_n n^{-3/2}$$

$$\leq \sqrt{c_4} \sum_{k=1}^{\infty} \frac{\sqrt{u_k}}{k^2} \sum_{n=u_k+1}^{u_{k+1}} n^{-3/2} + \sqrt{c_4} \sum_{k=1}^{\infty} \frac{1}{u_{k+1}} \frac{1}{k^2} \sum_{n=u_k+1}^{u_{k+1}} 1$$

$$+ \sqrt{c_4} \sum_{k=1}^{\infty} \frac{1}{\sqrt{u_{k+1}}} \frac{1}{k^2} \sum_{n=u_k+1}^{u_{k+1}} n^{-1/2} < \infty.$$

To prove that the second sum is finite, we need to analyze $A_n$, which satisfies the renewal equation

$$A_n = E_0[S_{n \wedge \nu}] + \sum_{j=1}^{n-1} A_{n-j} p_j.$$

Unlike Isola ([1999](#)), we use probabilistic arguments to analyze this renewal equation.

We define

$$T_0 = 0, \qquad T_k = \min\{t > T_{k-1} : Y_t = 0\}, \qquad \tau_k = T_k - T_{k-1}, \qquad k = 1, 2, \ldots.$$

Then, $\{\tau_j\}$ are independent variables equally distributed as $\tau$. [See, e.g., Breiman ([1968](#)), page 146.] Let $\xi_j = 1 - \pi_0 \tau_j$ and introduce the stopping time

$$\nu_n = \min\{j \geq 1 : T_j \geq n\}.$$

Clearly, $S_{T_k} = \sum_{j=1}^{k} \xi_j$, $E_0[\xi_1] = 0$, $\nu_n \leq n$ and, thus, by the Wald identity,

$$E_0[S_{T_{\nu_n}}] = E\left[\sum_{j=1}^{\nu_n} \xi_j\right] = 0.$$

Hence, since $|S_a - S_b| \leq |a - b|$, by the definition of $A_n$, we obtain

$$|A_n| = |E_0[S_{T_{\nu_n}} - S_n]| \leq E_0[\tau_{\nu_n}] \leq E_0\left[\max_{1 \leq i \leq n} \tau_i\right].$$

Let us denote by

$$M_n = \max_{1 \leq i \leq n} \tau_i.$$

Then,

$$J_n = \max_{1 \leq i \leq n} |A_i| \leq E[M_n].$$

To analyze $E[M_n]$, we notice that

$$E[M_n] = \sum_{t=1}^{\infty} u_t P(M_n = u_t)$$

and

$$P(M_n = u_t) \le \min(1, nP(\tau = u_t)) \le c_1 \min(1, n/u_t^2).$$

Fix $n$, $u_k < n \le u_{k+1}$. Notice first that, for $t \le k-1$, we have $u_t \le u_{k-1} \le u_k^{1/4} \le n^{1/4}$. Also, $\sum_{j=k+1}^{\infty} 1/u_j \le c_3/u_{k+1}$ and, thus, splitting the sum into three parts according to $t$: $t \le k-1$, $t = k$ and $t \ge k+1$, we obtain the bound

$$E[M_n] \le c_4 \left( n^{1/4} + \frac{n}{u_{k+1}} + u_k \min(1, n/u_k^2) \right).$$

Finally, by the construction of $u_n$ and its relation to $a_n$, we derive

$$\sum_{n=1}^{\infty} a_n J_n/n^{3/2} \le c_5 \sum_{n=1}^{\infty} n^{-5/4} + c_6 \sum_{k=1}^{\infty} \frac{1}{k^2} \frac{1}{u_{k+1}} \sum_{n=u_k+1}^{u_{k+1}} n^{-1/2}$$

$$+ c_7 \sum_{k=1}^{\infty} \frac{1}{k^2} \frac{1}{u_k} \sum_{n=u_k+1}^{u_k^2} n^{-1/2} + c_8 \sum_{k=1}^{\infty} \frac{1}{k^2} u_k \sum_{n=u_k^2+1}^{u_{k+1}} n^{-3/2}$$

$$< \infty,$$

proving (15).

3.3. *Stochastic unboundedness of $S_n/\sqrt{n}$ and the proof of Theorem 2.* We proceed by contradiction; that is, we assume that

$$\{S_n/\sqrt{n}; n \ge 1\} \qquad \text{is stochastically bounded}$$

and show that $E\tau^2 < \infty$, which is in contradiction with (14).

Let $\{\tau_j\}$ be independent variables equally distributed as $\tau$. Define

$$T_k = \tau_1 + \cdots + \tau_k, \qquad \eta_n = \max\{i \ge 1 : T_i \le n\},$$

$$T(i,n] = T_n - T_i, \qquad \eta_n(\xi) = \max\{i \ge 1 : \xi + T(1,i] \le n\}$$

(where $\max_{i \in \varnothing} a_i = 0$). Then, $S_n = \eta_n(\nu) - na$, where $a = 1/E[\tau_1] = \pi_0$.

The following proposition will provide a slightly more general result which has interest in itself.

PROPOSITION 3.2.  *Assume that, for a nonnegative integer valued variable $\xi$,*

(16)  $$\left\{ \frac{\eta_n(\xi) - an}{\sqrt{n}}; n \ge 1 \right\} \qquad \text{is stochastically bounded.}$$

*Then, $E[\tau_1^2] < \infty$.*

PROOF. First, let $\eta'_n$ be a copy of the renewal process $\{\eta_n : n \geq 1\}$ which does not depend on $\xi$. Then, $\eta_n(\xi)$ is equally distributed as $\eta'_{n-\xi}$ and so, any finite number of renewals do not affect the stochastic boundedness of the normalized renewal processes. As a consequence, condition (16) implies that

$$P([an - \sqrt{n}M] \leq \eta_n < [an + \sqrt{n}M]) \geq 1 - \varepsilon_M,$$

where $\varepsilon_M \to 0$ as $M \to \infty$.

Next, we apply the standard relationship $\{\eta_n \geq k\} = \{T_k \leq n\}$, yielding

$$P([an - \sqrt{n}M] \leq \eta_n < [an + \sqrt{n}M]) = P(T_{[an-\sqrt{n}M]} \leq n, T_{[an+\sqrt{n}M]} > n)$$
$$\equiv P(T_L \leq n, T_R > n) = I \geq 1 - \varepsilon_M,$$

where

$$L = L[n, M] = [an - \sqrt{n}M], \qquad R = R[n, M] = [an + \sqrt{n}M].$$

Now, we take $k = R - L$. Since $T(i, n) = T_n - T_i$ is equally distributed as $T_{n-i}$, we can write

$$I = P(T_L \leq n, T_L + T(L, R) > n)$$
$$= P(T_L \leq n - kN, T_L + T(L, R) > n)$$
$$\quad + P(n - kN < T_L \leq n, T_L + T(L, R) > n)$$
$$\leq P(T(L, R) > kN) + P(n - kN < T_L \leq n)$$
$$= P(T_k > kN) + P(n - kN < T_L \leq n).$$

By the law of the large numbers,

$$P(T_k/k > N) \leq \delta_N,$$

where $\delta_N \to 0$ as $N \to \infty$. Thus,

$$P(n - kN < T_L \leq n) \geq 1 - \varepsilon_M - \delta_N.$$

Since $2\sqrt{n}M - 1 \leq k \leq 2\sqrt{n}M + 1$, we derive

$$P(|T_L - n|/\sqrt{n} \leq (2M + 1)N) \geq 1 - \varepsilon_M - \delta_N.$$

Now we use the symmetrization argument. We consider an independent copy of $\{\tau_j\}$, namely, $\{\tau'_j\}$ and denote by $T'_k = \tau'_1 + \cdots + \tau'_k$, $T^s_k = T_k - T'_k$. Clearly,

$$P(|T^s_L|/\sqrt{n} \geq 2(2M + 1)N) \leq \varepsilon_M + \delta_N.$$

Here

$$\lim_{n \to \infty} L(n, M)/n = a.$$

By standards arguments involving an application of the Lévy maximal inequality for sums of symmetric independent random variables, we easily derive that the sequence $\{T_n^s/\sqrt{n}\}$ is stochastically bounded. By Theorem 3 in Esseen and Janson (1985), the fact that $\{T_n^s/\sqrt{n}\}$ is stochastically bounded implies $E(\tau_1 - \tau_1')^2 < \infty$. Thus, $E\tau_1^2 < \infty$.  □

PROOF OF THEOREM 1.2.   By combining Proposition 3.1 with the bound (15), we obtain the first part of (3). To prove the second part, we proceed by absurd and notice that if $\{S_n/\sqrt{n}\}$ is stochastically bounded, then by Proposition 3.2, $E[\tau^2] < \infty$, which is in contradiction with (14).

The proof of Theorem 1.2 is complete.   □

**Acknowledgments.**  The authors would like to express their gratitude to the referee whose suggestions led to significant simplification of the proofs, contributing in this way to the improvement of the presentation of this paper. We would also want to thank to Wei Biao Wu for showing us an interesting alternative way to improve on the constants in our former version of this paper. His comments led us to refine our proofs to further improve the constants in the maximal inequality.

DEPARTMENT OF MATHEMATICAL SCIENCES
UNIVERSITY OF CINCINNATI
P.O. BOX 210025
CINCINNATI, OHIO 45221-0025
USA
E-MAIL: peligrm@math.uc.edu

SCHOOL OF MATHEMATICAL SCIENCES
UNIVERSITY OF NOTTINGHAM
NOTTINGHAM, NG7 2RD
ENGLAND
E-MAIL: pmzsu@gwmail.nottingham.ac.uk